\documentclass[dvipdfmx]{article}
\usepackage{amsmath, amssymb, amsthm, delarray}

\usepackage[top =20mm, left=35mm,right=37mm]{geometry}
\usepackage{graphicx}
\usepackage{tikz}
\usetikzlibrary{intersections,calc,arrows.meta}

\newtheorem{thm}{Theorem}[section]
\newtheorem*{theorem*}{Theorem}
\newtheorem{lem}[thm]{Lemma}
\newtheorem{prop}[thm]{Proposition}

\theoremstyle{definition}

\newtheorem{rem}[thm]{Remark}

\newtheorem{asmp}[thm]{Assumption}

\numberwithin{equation}{section}
\theoremstyle{remark}

\begin{document}

\title{\textbf{An explicit lower bound for large gaps between some consecutive primes}}

\author{Keiju Sono}

\date{}
\allowdisplaybreaks

\maketitle 
\noindent
\begin{abstract}
Let $p_{n}$ denote the $n$th prime and for any fixed positive integer $k$ and $X\geq 2$, put
\[
G_{k}(X):=\max _{p _{n+k}\leq X} \min \{ p_{n+1}-p_{n}, \ldots , p_{n+k}-p_{n+k-1} \}.
\]
Ford, Maynard and Tao \cite{FMT} proved that there exists an effective absolute constant $c_{LG}>0$ such that
\[
G_{k}(X)\geq \frac{c_{LG}}{k^{2}}\frac{\log X \log \log X \log \log \log \log X}{\log \log \log X}
\]
holds for any sufficiently large $X$.  The main purpose of this paper is to clarify the numerical value of the constant $c_{LG}$ such that the above inequality holds.  We see that $c_{LG}$ is determined by several factors related to analytic number theory, for example, the ratio of integrals of functions in the multidimensional sieve of Maynard \cite{May}, the distribution of primes in arithmetic progressions to large moduli, and the coefficient of upper bound sieve of Selberg. We prove that the above inequality is valid at least for $c_{LG}\approx 2.0\times 10^{-17}$. 

\end{abstract}

%%%%%%%%%%%%%%%%%%%%%%%%%%%%%%%%%%%%%%%%%%%%%%%%%%%%%%%%%%%%%%%%%%%%%%%%%%%%%%%%%%%%%%%%%%%%%%%%%%%%%%%%%%%%
%%%%%%%%%%%%%%%%%%%%%%%%%%%%%%%%%%%%%%%%%%%%%%%%%%%%%%%%%%%%%%%%%%%%%%%%%%%%%%%%%%%%%%%%%%%%%%%%%%%%%%%%%%%%

\section{Introduction}
For a positive integer $k$ and a real number $X\geq 3$, put
\[
G_{k}(X):=\max _{p _{n+k}\leq X} \min \{ p_{n+1}-p_{n}, \ldots , p_{n+k}-p_{n+k-1} \}
\]
and in particular, put $G(X):=G_{1}(X)$, which denotes the largest prime gap below $X$.  The prime number theorem yields the number of primes below $x$ is asymptotically $x/\log x$, so the average gap of consecutive primes below $x$ is $(1+o(1))\log x$. Hence we get the trivial lower bound 
\[
G(X)\geq (1+o(1))\log X.
\]
Harald Cram\'{e}r \cite{C} constructed a simple probabilistic model of the set of prime numbers, and following his model, he conjectured that $G(X)\sim (\log X)^{2}$. (See also refinements of Cram\'{e}r's model by Granville \cite{Gr}, Firoozbakht \cite{R3} and Wolf \cite{Wo}.) In 1931, Westzynthius \cite{W} made the first quantitative improvement on the trivial bound and proved
\[
G(X)\gg \frac{\log X\log _{3}X}{\log _{4}X},
\]
where $\log _{n}$ denotes the $n$ times composition of the logarithmic function, i.e., $\log _{n}X:=(\underbrace{\log \circ \cdots \circ \log}_{n} )X$. In particular, it follows that $G(X)/\log X \to \infty$ as $X\to \infty$, so $G(X)$ can be arbitrarily large compared with the average gap. In 1934, Ricci \cite{Ri} slightly improved this and showed $G(X)\gg \log X \log _{3}X$. Erd\H{o}s  \cite{E} improved on Westzynthius' result and obtained
\[
G(X)\gg \frac{\log X \log _{2}X}{(\log _{3}X)^{2}}.
\]
Rankin \cite{R} made a further improvement and showed that
\begin{equation}
\label{R}
G(X)\geq (c+o(1))\frac{\log X \log _{2}X \log _{4}X}{(\log _{3}X)^{2}}
\end{equation}
holds with $c=\frac{1}{3}$. Several mathematicians improved on the value of the coefficient $c$ above (see  Sch\"{o}nhage \cite{S} with $c=\frac{1}{2}e^{\gamma}$, Rankin \cite{R2} with $c=e^{\gamma}$, Maier and Pomerance \cite{MP} with $c=1.31256e^{\gamma}$, and Pintz \cite{P} with $c=2e^{\gamma}$).  Erd\H{o}s conjectured that (\ref{R}) holds with arbitrarily large $c$. This problem had been unsolved for a long time, but Maynard \cite{May} and a team of Ford, Green, Konyagin, Tao \cite{FGKT}  independently solved at almost the same time in August  2014. Several months later, in a joint work \cite{FGKMT}, they  obtained
\begin{equation}
\label{longgap}
G(X)\gg \frac{\log X\log _{2}X\log _{4}X}{\log _{3}X}.
\end{equation}
This was the first quantitative improvement of Rankin's bound (\ref{R}) in almost 80 years. 

Let $k$ be a fixed positive integer. The $G_{k}(X)$ above has also been studied for a long time. Following the argument of Cram\'{e}r,  it is conjectured that $G_{k}(X)\asymp \frac{1}{k}\log ^{2}X$.     Erd\H{o}s \cite{E} considered the case $k=2$ and proved
\[
G_{2}(X)/\log X \to \infty
\]
as $X \to \infty$. Maier \cite{Mai} showed
\[
G_{k}(X _{n})\gg _{k}\frac{\log X _{n} \log _{2}X_{n} \log _{4}X_{n}}{(\log _{3}X_{n})^{2}}
\]
by combining  his famous Maier matrix method and Pintz's ideas in \cite{P}, where $(X_{n})$ is some monotonically  increasing sequence such that $X_{n}\to \infty$ as $n\to \infty$.   Pintz \cite{P2} improved on Maier's result and proved 
\[
G_{k}(X_{n})/ \left( \frac{\log X_{n}\log _{2}X_{n}\log _{4}X_{n}}{(\log _{3}X_{n})^{2}} \right) \to \infty .
\]
An issue of Maier's argument is that one has to avoid (possible) Siegel zeros, and this is the reason why the results of Maier and Pintz above are restricted to a special sequence $(X_{n})$, rather than all sufficiently large $X$. 

In \cite{FMT}, Ford, Maynard and Tao succeeded in handling this difficulty and gave a lower bound for $G_{k}(X)$ for any sufficiently large $X$. Concretely, they proved that

\[
G_{k}(X)\gg \frac{1}{k^{2}} \frac{\log X \log _{2}X \log _{4}X}{\log _{3}X}
\]
holds as $X\to \infty$, and the implied constant above is absolute and effective.
In other words, for any sufficiently large $X\geq 3$, there exists an absolute computable constant $c_{LG}>0$ such that 
\begin{equation}
\label{1}
G_{k}(X)\geq \frac{c_{LG}}{k^{2}} \frac{\log X \log _{2}X \log _{4}X}{\log _{3}X}
\end{equation}
holds. In addition to the technique to handle the possibility of the existence of exceptional zeros of $L$-functions, their result relies on hypergraph covering theorem, the construction of multidimensional sieve weight, and very clever and highly technical probabilistic arguments in \cite{FGKMT}.

Though the lower bound in their theorem is effective, the explicit value of $c_{LG}$ (resp. the implied constant of (\ref{longgap})) is not mentioned in \cite{FMT} (resp. \cite{FGKMT}). (In the blog on the paper  \cite{FGKMT}, Tao says “we manage to avoid the use of the (ineffective) Siegel-Walfisz theorem by deleting an exceptional prime from the multidimensional Selberg sieve, leading to an effective (but quite small) value of $c$.” )  However,  clarifying the value of  $c_{LG}$ seems to be  significant because it has the merit of making clear which parameter affects the coefficient $c_{LG}$ and how. With this reason, the main purpose of this paper is to clarify what value is appropriate as the coefficient $c_{LG}$ in (\ref{1}) and how several factors in analytic number theory are related.  It is fully expected that this attempt might be  a clue to establish a new quantitative improvement on the lower bound of $G_{k}(X)$.  We prove that one can take
\[
c_{LG}= \frac{C_{PAP}^{2}\theta c_{I,J}e^{4\gamma}}{737280000 \log 5 \;  C_{UB}M(1+D_{PAP}^{-1})^{4}(25C_{UB}+20e^{\gamma}M)}  
\]
(see Theorem 3.6), where $\theta$ is a parameter which describes how primes are equidistributed in arithmetic progressions (see \cite{May}, Hypothesis 1),  $c_{I,J}$ is the ratio of integrals of functions (see (\ref{5})) in the multidimensional sieve of Maynard \cite{May}, $C_{PAP}$ and $D_{PAP}$ are constants in the statement on distribution of primes in arithmetic progressions to large moduli (see Assumption 3.2), and $C_{UB}$, $D_{UB}$ are the constants of upper bound sieve of Selberg (see Assumption 3.4) and $M:=\max \{D_{PAP}, D_{UB} \}$.  In Sections 2-3, we will see that unconditionally one can take $\theta =\frac{1}{3}$, $c_{I,J}=\frac{1}{4}$, $C_{PAP}=1-e^{-2}$, $D_{PAP}=160$, $C_{UB}=8e^{2\gamma}$ and $D_{UB}>0$ arbitrarily. 
 Numerically, one can take
 \[
 c_{LG}\approx 2.0\times 10^{-17}.
 \]
Obviously, this value is not the best that can be obtained with current techniques and could be improved to some extent with some effort, for example, by applying recent results on zero density estimates and explicit zero free regions of Dirichlet $L$-functions. In any case, however, as long as we rely on current methods, the coefficient is likely to have to be fairly small.

%%%%%%%%%%%%%%%%%%%%%%%%%%%%%%%%%%%%%%%%%%%%%%%%%%%%%%%%%%%%%%%%%%%%%%%%%%%%%%%%%%%%%%%%%%%%%%%%%%%%%%%%%%%%
%%%%%%%%%%%%%%%%%%%%%%%%%%%%%%%%%%%%%%%%%%%%%%%%%%%%%%%%%%%%%%%%%%%%%%%%%%%%%%%%%%%%%%%%%%%%%%%%%%%%%%%%%%%%

\section{Notation}
Let $c>0$ be a fixed constant (to be determined later) and $x\geq 10$  a sufficiently large real number. Put
\begin{equation}
\label{2}
y:=c\frac{x\log x \log _{3}x}{\log _{2}x},   \quad z:=x^{\frac{\log _{3}x}{4\log _{2}x}}.
\end{equation}
Let $B_{0}$ be either $1$ or a prime number satisfying 
\begin{equation}
\label{3}
\log x \ll B_{0}\leq x.
\end{equation}
Define three disjoint sets of primes ${\cal S}$, ${\cal P}$ and ${\cal Q}$ by 
\[
{\cal S}:=\{s: \mathrm{prime} \; | \; \log ^{20}x <s \leq z, \; s \neq B_{0} \},
\]
\[
{\cal P}:=\{p: \mathrm{prime} \; | \; \frac{x}{2} <p \leq x, \; p \neq B_{0} \},
\]
\[
{\cal Q}:=\{q: \mathrm{prime} \; | \; x <q \leq y, \; q \neq B_{0} \}.
\]
For vectors of residue classes $\vec{a}=(a_{s}\; \textrm{mod}\; s)_{s\in {\cal S}}$, $\vec{n}=(n_{p}\; \textrm{mod}\; p)_{p\in {\cal P}}$, put 
\[
S(\vec{a}):=\{n \in \mathbb{Z}\; | \; n \not \equiv a_{s}\; (\mathrm{mod}\; s), \; \forall s \in {\cal S} \},
\]
\[
S(\vec{n}):=\{n \in \mathbb{Z}\; | \; n \not \equiv n_{p}\; (\mathrm{mod}\; p), \; \forall p \in {\cal P} \}.
\]
Let $0<\theta <1$ be a parameter of Hypothesis 1 of \cite{May}. In our situation, unconditionally one can take
\[
\theta =\frac{1}{3}
\]
(see \cite{FGKMT}, Section 8). 

Next, for $r\in \mathbb{Z}_{\geq 2}$, we denote by ${\cal F}_{r}$ the set of square-integrable symmetric functions $F: \mathbb{R}^{r}\to \mathbb{R}$ supported in 
\[
{\cal R}_{r}:=\{(x_{1}, \ldots ,x_{r})\in \mathbb{R}^{r}\; | \; x_{1}, \; \ldots ,x_{r}\geq 0, x_{1}+\ldots +x_{r}\leq 1 \}.
\]
For $F\in {\cal F}_{r}$, put 
\[
I_{r}(F):=\int _{0}^{\infty}\cdots \int _{0}^{\infty}F(t_{1}, \ldots ,t_{r})^{2}dt_{1}\cdots dt _{r},
\]
\[
J_{r}(F):=\int _{0}^{\infty}\cdots \int _{0}^{\infty}\left(\int _{0}^{\infty}F(t_{1}, \ldots ,t_{r})dt_{r} \right)^{2}dt_{1}\cdots dt_{r-1}.
\]
We suppose that for some positive constant $c_{I,J}$, the inequality
\begin{equation}
\label{5}
\sup _{F\in {\cal F}_{r}}\frac{J_{r}(F)}{I_{r}(F)}\geq c_{I,J}\frac{\log r}{r}\left(1+O\left( \frac{1}{\log r} \right) \right)
\end{equation}
holds as $r\to \infty$. Maynard  showed that one can take $c_{I,J}=\frac{1}{4}$ (see (8.27) of \cite{May}).

%%%%%%%%%%%%%%%%%%%%%%%%%%%%%%%%%%%%%%%%%%%%%%%%%%%%%%%%%%%%%%%%%%%%%%%%%%%%%%%%%%%%%%%%%%%%%%%%%%%%%%%%%%%%
%%%%%%%%%%%%%%%%%%%%%%%%%%%%%%%%%%%%%%%%%%%%%%%%%%%%%%%%%%%%%%%%%%%%%%%%%%%%%%%%%%%%%%%%%%%%%%%%%%%%%%%%%%%%

\section{The sieve of intervals and the explicit lower bound for $G_{k}(X)$}
\begin{prop}
Let $A\geq 1$ be an arbitrarily fixed constant and $x\geq 2$  sufficiently large real number. Let $y$ be a parameter defined by (\ref{2}) with 
\begin{equation}
\label{6}
c=\frac{\theta c_{I,J}}{12800\log 5}.
\end{equation}
Let $B_{0}$ be either $1$ or a prime number with $\log x \ll B_{0}\leq x$. Then, for any prime number $p\leq x$ with $p\neq B_{0}$, there exists a residue class $a_{p}\; (\textrm{mod}\; p)$ for which the set
\[
{\cal T}:=\{n\in [y]\backslash [x] \; | \; n\not \equiv a_{p}\; (\mathrm{mod}\; p), \; \forall p \leq x, \; p\neq B_{0} \}
\]
satisfies the following three conditions. \\
$\cdot$ (Upper bound) 
\begin{equation}
\label{7}
\# {\cal T}\leq 5A(1+o(1))\frac{x}{\log x}.
\end{equation}
$\cdot$ (Lower bound) 
\begin{equation}
\label{8}
\# {\cal T}\geq A(1+o(1))\frac{x}{\log x}.
\end{equation}
$\cdot$ (Upper bound in short intervals) For any fixed $0\leq \alpha <\beta \leq 1$, we have
\begin{equation}
\label{9}
\# ({\cal T}\cap (\alpha y, \beta y])\leq 5A(2|\beta -\alpha|+\varepsilon)(1+o(1))\frac{x}{\log x}.
\end{equation}
\end{prop}
We will prove this proposition in subsequent sections. We now give a lower bound for $G_{k}(X)$ as a consequence of this proposition. Let $Q\geq 100$. Then it is known that there exists a positive integer $B_{Q}$ which is either $1$ or a prime number with
\begin{equation}
\label{10}
B_{Q}\gg \log _{2}Q
\end{equation}
such that if a Dirichlet character $\chi$ with conductor less than $Q$ and coprime to $B_{Q}$ satisfies 
\[
L(\sigma +it, \chi)=0 \quad (\sigma ,t \in \mathbb{R}),
\]
then 
\begin{equation}
\label{11}
1-\sigma \geq \frac{c_{ZFR}}{\log (Q(1+|t|))}
\end{equation}
holds for some absolute constant $c_{ZFR}>0$.  We introduce the following two assumptions on distribution of primes in arithmetic progressions. 

\begin{asmp}[Assumption PAP]
Suppose that all $L$-functions associated to Dirichlet characters $\chi$ modulo $q$ do not have any zero in the region (\ref{11}) (with $Q$ replaced by $q$). Then, there exist  absolute constants $0<C_{PAP}\leq 1$ and  $D_{PAP}\geq 1$ for which for any positive integers $a, q$ with $(a,q)=1$ and $x\geq q^{D_{PAP}}$, one has 
\begin{equation}
\label{12}
\# \{p: \mathrm{prime}\; | \; p\leq x, \; p\equiv a \; (\mathrm{mod}\; q) \} \geq C_{PAP} (1+o(1))\frac{x}{\varphi (q)\log x}.
\end{equation}
\end{asmp}
\begin{rem}
Later we will show that the above assumption is valid unconditionally for $D_{PAP}=160$, $C_{PAP}=1-e^{-2}$.  
\end{rem}
For $X>1$, we denote by $[X]$ the set of integers in the interval $[1, X]$. 
\begin{asmp}[Assumption UB]
There exist absolute constants $C_{UB}\geq 1$, $D_{UB}>0$ for which the following holds. For any sufficiently large $x$ and a positive integer $B_{0}$ which is either $1$ or a prime number with
\begin{equation}
\label{13}
\log x \ll B_{0}\leq x,
\end{equation}
put $P:=P(x)/B_{0}$, where $P(x):=\prod _{p\leq }p$ denotes the product of all primes equal to or less than $x$. Then, for any $Z\geq P^{D_{UB}}$ and $a, b \in [P]$ with $a\neq b$, one has 
\begin{equation}
\label{14}
\begin{aligned}
& \# \{z\in [Z] \; | \; Pz+a, \; Pz+b : \mathrm{prime} \} \leq C_{UB}(1+o(1)) \left( \frac{\log x}{\log Z} \right)^{2}Z.
\end{aligned}
\end{equation}
\end{asmp}

\begin{rem}
Later we will show that one can unconditionally take $C_{UB}=8e^{2\gamma}$ and $D_{UB}>0$ arbitrarily. 
\end{rem}

With these notations, we have the following theorem.
\begin{thm}[An explicit version of the theorem of Ford, Maynard and Tao \cite{FMT}]
Let $k$ be a fixed positive integer. Let $\theta$ be a parameter of Maynard's sieve in Hypothesis 1 of \cite{May}, $c_{I,J}$ a constant satisfying (\ref{5}), $C_{PAP}$ and $D_{PAP}$  constants in Assumption PAP and $C_{UB}, D_{UB}$ constants in Assumption UB and put $M:=\max \{D_{PAP}, D_{UB} \}$. Then we have
\begin{equation}
\label{MT}
\begin{aligned}
G_{k}(X)\geq & \frac{C_{PAP}^{2} \theta c_{I,J}e^{4\gamma}}{737280000 \log 5 \;  C_{UB}M(1+D_{PAP}^{-1})^{4}(25C_{UB}+20e^{\gamma}M)}  \frac{1}{k^{2}}\frac{\log X \log _{2}X\log _{4}X}{\log _{3}X}
\end{aligned}
\end{equation}
fot any sufficiently large $X$. 
\end{thm}

\begin{lem}
Assume Assumptions PAP and UB. Let $B_{0}$ be a positive integer which is either $1$ or a prime number satisfying (\ref{13}). Put $P=P(x)/B_{0}$. Then, for any $Z\geq P^{D_{PAP}}$ and $a\in [P]$ with $(a,P)=1$, we have
\begin{equation}
\label{15}
\# \{z\in [Z]\; | \; Pz+a: \mathrm{prime} \} \geq \frac{e^{\gamma}}{1+D_{PAP}^{-1}}(1+o(1))\frac{\log x}{\log Z}Z.
\end{equation}
\end{lem}

\begin{proof}
We apply the Assumption PAP with a zero free region (\ref{11}) with $Q=P(x)=\prod _{p\leq x}p$. In the Assumption UB, we set $B_{0}=1$ if $B_{P(x)}>x$ and otherwise put $B_{0}=B_{P(x)}$. Note that  in the latter case, due to the condition (\ref{10}) and a consequence of the prime number theorem ($P(x)\sim e^{(1+o(1))x}$), it follows that
\[
B_{0}=B_{P(x)}\gg \log _{2}P(x)=\log _{2}\prod _{p\leq x}p \sim \log _{2}e^{x}=\log x,
\]
so (\ref{13}) is satisfied. Therefore, by the condition (\ref{12}), for $Z\geq P^{D_{PAP}}$, we have
\begin{equation}
\label{16}
\begin{aligned}
& \# \{z\in [Z]\; | \; Pz+a: \mathrm{prime} \}  \\
&=\# \{P+a \leq p \leq PZ+a \; | \; p\equiv a\; (\mathrm{mod}\; P), \; p:\mathrm{prime} \} \\
&\geq C_{PAP}(1+o(1))\frac{PZ+a}{\varphi (P)\log (PZ+a)}-\# \{p<P+a \; | \; p\equiv a\; (\mathrm{mod}\; P), \; p:\mathrm{prime} \}.
\end{aligned}
\end{equation}
By Mertens' formula, we have
\begin{equation}
\label{Mer}
\frac{P}{\varphi (P)}=\frac{P}{P\prod _{p\leq x, p\neq B_{0}}\left(1-\frac{1}{p} \right)}=\frac{1}{\prod _{p\leq x, p\neq B_{0}}\left(1-\frac{1}{p} \right)}\sim e^{\gamma}\log x.
\end{equation}
Hence the first term of the last line of (\ref{16}) is asymptotically
\[
\frac{PZ}{\varphi (P)\log PZ}\geq e^{\gamma}\log x\frac{Z}{\log Z^{1+\frac{1}{D_{PAP}}}}=\frac{e^{\gamma}}{1+D_{PAP}^{-1}}\frac{\log x}{\log Z}Z.
\]
On the other hand, the second term of the last line of (\ref{16}) is at most
\[
\# \{p<P \; | \; p:\mathrm{prime} \} \leq \frac{2P}{\log P}\leq 2D_{PAP}\frac{Z^{\frac{1}{D_{PAP}}}}{\log Z}=o\left(\frac{\log x}{\log Z}Z \right).
\]
Combining these estimates, we obtain (\ref{15}). 
\end{proof}
We write $\mathbb{P}$ for probability and $\mathbb{E}$ for expectation. Put  $Z=P^{M}$, where $M:=\max \{ D_{UB}, D_{PAP} \}$.  Let $\mathbf{z}$ be a random variable in $[Z]$ which is chosen uniformly and $y$, ${\cal T}$ and $a_{p} \; (\textrm{mod}\; p)$ are those in Proposition 3.1. By Chinese remainder theorem, there exists $m \in [Z]$ such that $m \equiv -a_{p}\; (\textrm{mod}\; p)$ holds for any $p\leq x, p \neq B_{0}$.  Therefore, the interval $\mathbf{z}P+m+{\cal T}$ exactly  consists of elements of $\mathbf{z}P+m+[y]\backslash [x]$ which is coprime to $P$. In particular, all prime numbers in $\mathbf{z}P+m+[y]\backslash [x]$ are contained in $\mathbf{z}P+m+{\cal T}$. By (\ref{15}) and $Z=P^{M}=e^{(1+o(1))Mx}$, we have
\begin{equation}
\label{17}
\begin{aligned}
\mathbb{P}(\mathbf{z}P+m+a: \mathrm{prime})&\geq \frac{C_{PAP}e^{\gamma}}{1+D_{PAP}^{-1}}(1+o(1))\frac{\log x}{\log Z} \\
&=\frac{C_{PAP}e^{\gamma}}{(1+D_{PAP}^{-1})M}(1+o(1))\frac{\log x}{x}
\end{aligned}
\end{equation}
for any $a\in {\cal T}$. On the other hand, by Brun-Titchmarsh type estimate 
\[
\pi (x; q,a)\leq \frac{2x}{\varphi (q)\log \frac{x}{q}} \quad (\forall x >q)
\]
and Mertens type estimate (\ref{Mer}), we have
\begin{equation}
\label{18}
\begin{aligned}
\mathbb{P}(\mathbf{z}P+m+a: \mathrm{prime})&=\frac{1}{Z}\# \{z\in [Z] \; | \; Pz+m+a: \mathrm{prime} \} \\
&\leq \frac{1}{Z}\# \{p\leq PZ+m+a \; | \; p\equiv m+a \; (\mathrm{mod}\; P) \} \\
&\leq \frac{1}{Z}\# \{p\leq 2PZ \; | \; p\equiv m+a \; (\mathrm{mod}\; P) \} \\
&\leq \frac{1}{Z}\frac{2\cdot 2PZ}{\varphi (P)\log \frac{2PZ}{P}} \\
&\leq \frac{4e^{\gamma}\log x}{\log (2Z)} \\
&\leq \frac{4e^{\gamma}\log x}{\log Z}\sim \frac{4e^{\gamma}}{M}\frac{\log x}{x}.
\end{aligned}
\end{equation}
Furthermore, by Assumption UB, 
\begin{equation}
\label{19}
\begin{aligned}
\mathbb{P}(\mathbf{z}P+m+a, \; \mathbf{z}P+m+b: \mathrm{prime})&\leq C_{UB}(1+o(1)) \left( \frac{\log x}{\log Z} \right)^{2} \\
&=\frac{C_{UB}}{M^{2}}(1+o(1))\left( \frac{\log x}{x} \right)^{2}
\end{aligned}
\end{equation}
for any $a\neq b\in {\cal T}$. Let $\mathbf{N}$ denote the number of primes in $\mathbf{z}P+m+{\cal T}$ (hence in $\mathbf{z}P+m+[y]\backslash [x]$). By (\ref{17}) and (\ref{8}), we have
\begin{equation}
\label{20}
\begin{aligned}
\mathbb{E}\mathbf{N}&=\sum _{a\in {\cal T}}\mathbb{P}(\mathbf{z}P+m+a: \mathrm{prime}) \geq (1+o(1))\frac{C_{PAP}e^{\gamma}}{(1+D_{PAP}^{-1})M}A.
\end{aligned}
\end{equation}
On the other hand, by (\ref{18}), (\ref{19}) and (\ref{7}), 
\begin{equation}
\label{21}
\begin{aligned}
\mathbb{E}\mathbf{N}^{2}&=\sum _{a, b\in {\cal T}}\mathbb{P}(\mathbf{z}P+m+a, \mathbf{z}P+m+b: \mathrm{prime}) \\
&\leq \frac{C_{UB}}{M^{2}}(1+o(1))\left( \frac{\log x}{x} \right)^{2}(\# {\cal T})^{2}+\frac{4e^{\gamma}}{M}(1+o(1))\frac{\log x}{x}\# {\cal T} \\
&\leq \left( \frac{25A^{2}C_{UB}}{M^{2}}+\frac{20Ae^{\gamma}}{M} \right)(1+o(1)) \\
&\leq \frac{25C_{UB}+20e^{\gamma}M}{M^{2}}(1+o(1))A^{2}.
\end{aligned}
\end{equation}
(The first term of the second line is  an upper bound for the contribution of the terms with $a\neq b$, and the second term is that of the terms with $a=b$. In the last line, we used $A\geq 1$.) Suppose that $0\leq \alpha <\beta \leq 1$ satisfy $\beta -\alpha \leq 2\varepsilon$ for $0<\varepsilon <1$. Then by (\ref{19}) and (\ref{9}), the probability that the interval $\mathbf{z}P+m+[\alpha y, \beta y]$ contains at least two primes is at most
\[
\left(6A(4\varepsilon +\varepsilon )\frac{x}{\log x} \right)^{2}\frac{C_{UB}}{M^{2}}(1+o(1))\left( \frac{\log x}{x} \right)^{2}\leq \frac{900\varepsilon ^{2}A^{2}C_{UB}}{M^{2}}(1+o(1)).
\]
We cover the interval $[0, 1]$ by intervals
\[
I_{1}=[0, 2\varepsilon ], \;   I_{2}=[\varepsilon , 3\varepsilon ], \;  \ldots , \; I_{n_{\varepsilon}}=[(n_{\varepsilon}-1)\varepsilon ,1],
\]
where $n_{\varepsilon}=\lceil \frac{1}{\varepsilon} \rceil $. Then two elements $a, b \in [0, 1]$ with $|a-b|\leq \varepsilon$ are contained in the same interval $I_{j}$ for some $1\leq j \leq n_{\varepsilon}$. Therefore, the probability that “The interval $\mathbf{z}P+m+[y]\backslash [x]$ contains two primes with gap at most $\varepsilon y$ ”  is equal to or less than the probability that “some of the $n_{\varepsilon}$ intervals
\[
\mathbf{z}P+m+[0, 2\varepsilon y], \; \mathbf{z}P+m+[\varepsilon y, 3\varepsilon y], \; \ldots , \; \mathbf{z}P+m+[(n_{\varepsilon}-1)\varepsilon y, y]
\]
contains two primes”,  and by the above consideration, this probability is at most
\[
n_{\varepsilon}\frac{900\varepsilon ^{2}A^{2}C_{UB}}{M^{2}}(1+o(1)) <\frac{1800 \varepsilon A^{2}C_{UB}}{M^{2}},
\]
since $n_{\varepsilon}=\left\lceil \frac{1}{\varepsilon} \right\rceil <\frac{2}{\varepsilon}$. We take
\begin{equation}
\label{varepsilon}
\varepsilon =\frac{M^{2}}{1800A^{2}C_{UB}}\frac{e^{2\gamma}}{4(1+D_{PAP}^{-1})^{2}(25C_{UB}+20e^{\gamma}M)}.
\end{equation}
(Later we will see that $\varepsilon <1$ with an appropriate choice of $A$.)  Then the above probability is less than 
\[
\frac{e^{2\gamma}}{4(1+D_{PAP}^{-1})^{2}(25C_{UB}+20e^{\gamma}M)}.
\]
In other words, all primes in $\mathbf{z}P+m+[y]\backslash [x]$ are separated more than $\varepsilon y$ to each other with probability greater than
\begin{equation}
\label{22}
1-\frac{e^{2\gamma}}{4(1+D_{PAP}^{-1})^{2}(25C_{UB}+20e^{\gamma}M)}.
\end{equation}
On the other hand, for any fixed $c_{1}>0$, by Cauchy-Schwarz inequality, we have
\begin{align*}
\mathbb{E}\mathbf{N}&\leq c_{1}A\mathbb{P}(\mathbf{N}\leq c_{1}A)+\mathbb{E}\mathbf{N}1_{\mathbf{N}>c_{1}A} \leq c_{1}A+(\mathbb{E}\mathbf{N}^{2})^{\frac{1}{2}}\mathbb{P}(\mathbf{N}>c_{1}A)^{\frac{1}{2}}. 
\end{align*}
Therefore, by (\ref{20}) and (\ref{21}), we have
\begin{equation}
\label{Tao}
\begin{aligned}
\mathbb{P}(\mathbf{N}>c_{1}A)&\geq \frac{(\mathbb{E}\mathbf{N}-c_{1}A)^{2}}{\mathbb{E}\mathbf{N}^{2}} \geq \frac{\left(\frac{C_{PAP}e^{\gamma}A}{(1+D_{PAP}^{-1})M}-c_{1}A \right)^{2}}{\frac{(25C_{UB}+20e^{\gamma}M)A^{2}}{M^{2}}} =\frac{\left(\frac{C_{PAP}e^{\gamma}}{1+D_{PAP}^{-1}}-c_{1}M \right)^{2}}{25C_{UB}+20e^{\gamma}M},
\end{aligned}
\end{equation}
provided that $c_{1}<\frac{C_{PAP}e^{\gamma}}{(1+D_{PAP}^{-1})M}$. We take
\[
c_{1}=\frac{C_{PAP}e^{\gamma}}{2(1+D_{PAP}^{-1})M}.
\]
Then, 
\begin{equation}
\label{23}
\mathbb{P}(\mathbf{N}>c_{1}A)\geq \frac{e^{2\gamma}}{4(1+D_{PAP}^{-1})^{2}(25C_{UB}+20e^{\gamma}M)}.
\end{equation}
Since the sum of (\ref{22}) and the right hand side of (\ref{23}) is equal to $1$, it follows that there exists some integer $z\in [Z]$ such that the interval $zP+m+[y]\backslash [x]$ contains at least
\[
c_{1}A=\frac{e^{\gamma}A}{2M(1+D_{PAP}^{-1})}
\]
primes and all of them are separated at least
\[
\varepsilon y =\frac{yM^{2}}{1800A^{2}C_{UB}}\frac{e^{2\gamma}}{4(1+D_{PAP}^{-1})^{2}(25C_{UB}+20e^{\gamma}M)}
\]
to each other, where 
\[
y=c\frac{x\log x \log _{3}x}{\log _{2}x}, \quad c=\frac{\theta c_{I,J}}{12800\log 5}.
\]
We take
\begin{equation}
\label{24}
A=2e^{-\gamma}C_{PAP}^{-1}M(1+D_{PAP}^{-1})k.
\end{equation}
Then, the interval $zP+m+[y]\backslash [x]$ contains $k$ consecutive primes which are separated at least $\varepsilon y$ to each other. Therefore,
\begin{align*}
G_{k}(ZP+m+y)&\geq \varepsilon y \\
&=\frac{\theta c_{I,J}}{12800 \log 5}\frac{M^{2}}{1800A^{2}C_{UB}}\frac{e^{2\gamma}}{4(1+D_{PAP}^{-1})^{2}(25C_{UB}+20e^{\gamma}M)}\frac{x\log x\log _{3}x}{\log _{2}x}.
\end{align*}
Since
\[
ZP+m+y \leq 2ZP=2P^{M+1}=2e^{(1+o(1))(M+1)x}\leq e^{1.5Mx},
\]
by putting $X=e^{1.5Mx}$, we have
\begin{align*}
G_{k}(X)& \geq \frac{\theta c_{I,J}}{12800 \log 5}\frac{M^{2}}{1800A^{2}C_{UB}}\frac{e^{2\gamma}(1+o(1))}{4(1+D_{PAP}^{-1})^{2}(25C_{UB}+20e^{\gamma}M)} \frac{2}{3M}\frac{\log X \log _{2}X \log _{4}X}{\log _{3}X} \\
& \geq  \frac{\theta c_{I,J}}{12800 \log 5}\frac{M^{2}}{1800A^{2}C_{UB}}\frac{e^{2\gamma}}{4(1+D_{PAP}^{-1})^{2}(25C_{UB}+20e^{\gamma}M)} \frac{1}{2M} \frac{\log X \log _{2}X \log _{4}X}{\log _{3}X}.
\end{align*}
Finally, by substituting (\ref{24}), we obtain (\ref{MT}).   \hspace{6.5cm} $\Box$

%%%%%%%%%%%%%%%%%%%%%%%%%%%%%%%%%%%%%%%%%%%%%%%%%%%%%%%%%%%%%%%%%%%%%%%%%%%%%%%%%%%%%%%%%%%%%%%%%%%%%%%%%%%%
%%%%%%%%%%%%%%%%%%%%%%%%%%%%%%%%%%%%%%%%%%%%%%%%%%%%%%%%%%%%%%%%%%%%%%%%%%%%%%%%%%%%%%%%%%%%%%%%%%%%%%%%%%%%

\section{Possible values of $C_{UB}$ and $D_{UB}$ in Assumption UB}
To establish the upper bound (\ref{14}) in Assumption UB, we use the following theorem.
\begin{thm}[\cite{HR}, Theorem 5.7]
For $g\in \mathbb{N}$ and integers $a_{i}, b_{i}$ $(i=1, \ldots ,g)$, suppose
\[
E:=\prod _{i=1}^{g}a_{i}\prod _{1\leq r<s\leq g}(a_{r}b_{s}-a_{s}b_{r})\neq 0.
\]
For any prime number $p$, let $\rho (p)$ denote the number of solutions $n \; (\textrm{mod}\; p)$ of the equation 
\[
\prod _{i=1}^{g}(a_{i}n+b_{i}) \equiv 0 \quad (\mathrm{mod}\; p),
\]
and suppose that $\rho (p)<p$ holds for any primes $p$. Then, for any $1<y\leq x$, we have
\begin{align*}
& \#\{n \; | \; x-y<n \leq x, a_{i}n+b_{i}: \mathrm{prime} \} \\
& \quad \leq 2^{g}g!\prod _{p}\left(1-\frac{\rho (p)-1}{p-1} \right)\left(1-\frac{1}{p} \right)^{-g+1}\frac{y}{\log ^{g}y}  \left(1+O((\log y)^{-1}(\log _{2}3y+\log _{2}3|E|))   \right).
\end{align*}
\end{thm}
We apply this theorem with $x=y=Z$, $a_{1}=a_{2}=P$, and $b_{1}=a$, $b_{2}=b$.  Then, $E=P^{3}(a-b)\neq 0$ and 
\[
|E|\ll P^{3}Z\ll e^{3(1+o(1))x}Z.
\]
Note that since $a, b \in {\cal T}\subset [y]\backslash [x]$, we have $|a-b|\leq y =o(x\log x)$. We will use this later (see (\ref{30}) below). 

We need to compute $\rho (p)$. First, for $p\leq x$ with $p\neq B_{0}$, if either $p|a$ or $p|b$ holds, then
\[
Pn+a \equiv 0 \quad (\mathrm{mod}\; p) \quad \mathrm{or} \quad Pn+b \equiv 0 \quad (\mathrm{mod}\; p) 
\]
holds for any $n\in \mathbb{N}$, so 
\[
\# \{z \in [Z] \; | \; Pz+a, Pz+b :\mathrm{prime} \}=0.
\]
Therefore, it suffices to assume that $p\nmid ab$ ($\forall p\leq x$) or $p=B_{0}$ (in case of $B_{0}$: prime).  \\
I) If $p\leq x$, $p\neq B_{0}$, then the condition 
\[
(Pn+a)(Pn+b)\equiv 0 \quad (\mathrm{mod}\: p)
\]
is equivalent to 
\[
ab \equiv 0 \quad (\mathrm{mod}\: p)
\]
which does not hold due to our assumption. Hence $\rho (p)=0$. \\
II) If $B_{0}$ is a prime number and $p=B_{0}$, then the equation 
\[
(Pn +a)(Pn+b)\equiv 0 \quad (\mathrm{mod}\; p)
\]
has two solutions. Hence $\rho (B_{0})=2$. \\
III) If $p>x$, then the number of the solutions $n\; (\textrm{mod}\; p)$ of the equation
\[
(Pn +a)(Pn+b)\equiv 0 \quad (\mathrm{mod}\; p)
\]
is $1$ if $a\equiv b\; (\textrm{mod}\; p)$, and 2 if $a \not \equiv b\; (\textrm{mod}\; p)$. \\
Consequently, 
\begin{equation}
\label{26}
\begin{aligned}
& \# \{z\in [Z] \; | \; Pz+a, Pz+b :\mathrm{prime} \} \\
& \quad \leq 2^{2}\cdot 2! \underset{p\neq B_{0}}{\prod _{p\leq x}}\left(1+\frac{1}{p-1} \right)\left(1-\frac{1}{p} \right)^{-1}\left[\left(1-\frac{1}{B_{0}-1} \right)\left(1-\frac{1}{B_{0}} \right)^{-1}   \right] \\
& \quad \quad \quad  \times \underset{p\nmid (a-b)}{\prod _{p>x}}\left(1-\frac{1}{p-1} \right)\left(1-\frac{1}{p} \right)^{-1}\underset{p|(a-b)}{\prod _{p>x}}\left(1-\frac{1}{p} \right)^{-1} \\
& \quad \quad  \times \frac{Z}{\log ^{2}Z}\left(1+O\left( \frac{\log _{2}Z+\log x}{\log Z} \right) \right).
\end{aligned}
\end{equation}
Here, the factor $[\cdots]$ in the second line  appears only when $B_{0}$ is prime. Now, 
\begin{equation}
\label{27}
\underset{p\neq B_{0}}{\prod _{p\leq x}}\left(1+\frac{1}{p-1} \right)\left(1-\frac{1}{p} \right)^{-1} \leq \prod _{p\leq x}\left(1-\frac{1}{p} \right)^{-2}\sim e^{2\gamma}\log ^{2}x,
\end{equation}
\begin{equation}
\label{28}
\left(1-\frac{1}{B_{0}-1} \right)\left(1-\frac{1}{B_{0}} \right)^{-1}\sim 1, 
\end{equation}
\begin{equation}
\label{29}
\underset{p\nmid (a-b)}{\prod _{p>x}}\left(1-\frac{1}{p-1} \right)\left(1-\frac{1}{p} \right)^{-1}\leq 1, 
\end{equation}
and 
\begin{equation}
\label{30}
\underset{p| (a-b)}{\prod _{p>x}} \left(1-\frac{1}{p} \right)^{-1}=\exp \left(-\underset{p| (a-b)}{\sum _{p>x}}\log \left(1-\frac{1}{p}\right) \right) \leq \exp \left (\underset{p| (a-b)}{\sum _{p>x}} \frac{2}{p}\right) \leq \exp \left(\frac{2\omega (|a-b|)}{x}\right) \sim 1,
\end{equation}
where $\omega (n)$ denotes the number of distinct prime factors of $n$. Substituting  (\ref{27})-(\ref{30}) into (\ref{26}), we have
\begin{equation}
\begin{aligned}
\# \{z\in [Z] \; | \; Pz+a, Pz+b :\mathrm{prime} \} &\leq 8e^{2\gamma}\log ^{2}x \frac{Z}{\log ^{2}Z}(1+o(1))  \\
& =8e^{2\gamma}(1+o(1)) \left( \frac{\log x}{\log Z} \right)^{2}Z,
\end{aligned}
\end{equation}
provided that $\log x=o(\log Z)$. Since $Z\geq P^{D_{UB}}=e^{(1+o(1))D_{UB}x}$, this condition is valid for any $D_{UB}>0$. Consequently, the Assumption UB is valid for 
\begin{equation}
\label{30.5}
C_{UB}=8e^{2\gamma}
\end{equation}
and arbitrary $D_{UB}>0$.

%%%%%%%%%%%%%%%%%%%%%%%%%%%%%%%%%%%%%%%%%%%%%%%%%%%%%%%%%%%%%%%%%%%%%%%%%%%%%%%%%%%%%%%%%%%%%%%%%%%%%%%%%%%%
%%%%%%%%%%%%%%%%%%%%%%%%%%%%%%%%%%%%%%%%%%%%%%%%%%%%%%%%%%%%%%%%%%%%%%%%%%%%%%%%%%%%%%%%%%%%%%%%%%%%%%%%%%%%

\section{Possible values of $C_{PAP}$ and  $D_{PAP}$ in Assumption PAP}
To establish an inequality in Assumption PAP, we need some results on zero-free regions of Dirichlet $L$-functions. We use the following results by McCurley \cite{Mcc}. For a positive integer $q$, let ${\cal L}_{q}(s)$ be the product of $\varphi (q)$ $L$-functions associated to Dirichlet characters modulo $q$. Write $s=\sigma +it$, where $\sigma ,t\in \mathbb{R}$. 
\begin{prop}[\cite{Mcc}, Theorem 1]
Let $M:=\max \{q, q|t|, 10 \}$, $R=9.6459\ldots$.  Then ${\cal L}_{q}(s)$ has at most one zero in 
\[
\left\{s=\sigma +it \; | \; \sigma \geq 1-\frac{1}{R\log M} \right\}.
\]
The only possible zero is  a real zero of the $L$-function associated to a non-principal real character modulo $q$. 
\end{prop}
\begin{prop}[\cite{Mcc}, Theorem 2]
For $i=1,2$, let $\chi _{i}$ be a  real character modulo $q_{i}$ such that $\chi _{1}\neq \chi _{2}$ and $\beta _{i}$ be a real zero of $L(s, \chi _{i})$. Put
\[
M_{1}=\max \left\{ \frac{q_{1}q_{2}}{17}, 13 \right\}, \quad R_{1}=\frac{5-\sqrt{5}}{15-10\sqrt{2}}.
\]
Then we have
\[
\min \{\beta _{1}, \beta _{2} \}<1-\frac{1}{R_{1}\log M_{1}}.
\]
\end{prop}
For $Q\geq 100$ and $q_{1}, q_{2}\leq Q$, let $\chi _{1}$, $\chi _{2}$ be distinct primitive characters modulo $q_{1}$, $q_{2}$ respectively. Suppose that for $j=1,2$, $L(s, \chi _{j})$ has a zero $\rho _{j}=\beta _{j}+i\gamma _{j}$ in the domain
\[
{\cal R}: \left\{s=\sigma +it \; | \; \sigma \geq 1-\frac{1}{4R_{1}\log (Q(1+|t|))} \right\}.
\]
Then for $M_{(j)}:=\max \{q_{j}, q_{j}|t|, 10 \}$ ($j=1,2$), 
\[
1-\frac{1}{R\log M_{(j)}}\leq 1-\frac{1}{4R_{1}\log (Q(1+|t|))}
\]
holds. Hence by Proposition 5.1, $\rho _{j}$ is a real zero of $L(s, \chi_{j})$ for $j=1,2$. By Proposition 5.2, we have
\begin{equation}
\label{betabeta}
\min \{\beta _{1}, \beta _{2} \}\leq 1-\frac{1}{R_{1}\log M_{1}},
\end{equation}
where $M_{1}=\max \{\frac{q_{1}q_{2}}{17}, 13 \}$. However, since
\[
1-\frac{1}{R_{1}\log M_{1}}\leq 1-\frac{1}{4R_{1}\log Q},
\]
the inequality (\ref{betabeta}) contradicts our assumption that both $\rho _{1}=\beta _{1}$ and  $\rho _{2}=\beta _{2}$ lie in ${\cal R}$. Combining this and an argument in the proof of Corollary 1 of \cite{FMT}, we have the following conclusion. 
\begin{prop}
Put
\[
c_{ZFR}:=\frac{1}{24} \left(<\frac{1}{4R_{1}}=\frac{15-10\sqrt{2}}{4(5-\sqrt{5})} \right).
\]
Let $Q\geq 100$ be a fixed positive integer. Suppose that $L(s, \chi)=0$ holds for some primitive Dirichlet character  modulo at most $Q$. Then either  one of
\[
1-\sigma \geq \frac{c_{ZFR}}{\log (Q(1+|t|))} \quad  (s=\sigma +it, \sigma ,t \in \mathbb{R})
\]
or “$t=0$ and $\chi =\chi _{Q}$ is a real character which is determined uniquely by $Q$” holds. 
\end{prop}
Gallagher \cite{G} proved that there exist some absolute constants $a, b>0$ for which 
\begin{equation}
\label{G1}
\sum _{q\leq Q}\sum _{\chi (\mathrm{mod}\; q)}^{\quad \quad *}\left| \sum _{p=x}^{x+h}\chi (p)\log p \right| \ll h\exp \left(-a\frac{\log x}{\log Q} \right)
\end{equation}
holds uniformly for $x/Q \leq h \leq x$, $\exp (x^{\frac{1}{2}}\log x)\leq Q \leq x^{b}$. We need to make clear the value of $a$. Due to the above consideration, the $L$-functions associated to primitive characters modulo at most $T$ do not have zero in 
\[
\sigma >1-\frac{3c_{ZFR}}{\log T}, \quad |t|\leq T
\]
with at most one exception, and if the exceptional zero exists, then it is a simple real zero. Hence one can take $c_{1}=3c_{ZFR}$ in \cite{G}, (26). The estimate of \cite{G}, (30) yields
\begin{equation}
\label{G2}
\sum _{q\leq Q}\sum _{\chi (\mathrm{mod}\; q)}^{\quad \quad *}\left| \sum _{p=x}^{x+h}\chi (p)\log p \right| \ll h\left(\sum _{q\leq Q}\sum _{\chi (\mathrm{mod}\; q)}^{\quad \quad *}\sum _{\beta}x^{\beta -1}+\frac{Q^{4}}{T} \right),
\end{equation}
where $\beta$ runs over zeros of non trivial zeros of $L(s, \chi)$ in $\sigma \leq 1-3c_{ZFR}/\log T$, $|t|\leq T$. Due to the argument of \cite{G}, we have
\[
\sum _{q\leq Q}\sum _{\chi (\mathrm{mod}\; q)}^{\quad \quad *}\sum _{\beta}x^{\beta -1}\ll x^{-\frac{1}{2}\frac{3c_{ZFR}}{\log T}}.
\]
Therefore, by putting $T=Q^{5}$, the right hand side of (\ref{G2}) is at most
\[
\exp \left(-\frac{3}{10}c_{ZFR}\frac{\log x}{\log Q} \right).
\]
Therefore, one can take
\[
a=\frac{3}{10}c_{ZFR}
\]
in (\ref{G1}).

Next, by  Theorem 6 of \cite{G}, for some constant $c_{ZD}>0$, 
\begin{equation}
\label{31}
\sum _{q\leq T}\sum _{\chi \; (\mathrm{mod}\; q)}^{\quad \quad *}N(\alpha ,T, \chi)\ll T^{c_{ZD}(1-\alpha)}
\end{equation}
holds uniformly for $0\leq \alpha \leq 1$. Here, $N(\alpha ,T, \chi)$ denotes the number of zeros of $L(s, \chi)$ in the domain $\alpha \leq \Re (s)\leq 1$, $|\Im (s)|\leq T$. With this $c_{ZD}$, Theorem 7 of \cite{G} is valid for $T^{c_{ZD}}\leq x^{\frac{1}{2}}$, $T=Q^{5}$. Therefore, for $\exp (\log ^{\frac{1}{2}}x)\leq Q \leq x^{\frac{1}{10c_{ZD}}}$, (\ref{G1}) holds for $x/Q \leq h \leq x$. Therefore, the proof og Lemma 2 of Maier \cite{Mai} is valid for $x\geq Q^{10c_{ZD}}$.

In \cite{Mai}, the contribution of the terms with $\chi =\chi _{0}$ to $\psi (x; q,a)$ is $(1+o(1))x/\varphi (q)$, whereas that of the terms with $\chi \neq \chi _{0}$ is $O_{\leq}(xe^{-aD}/\varphi (q))$. Here, $f(x)=O_{\leq}(g(x))$ means that $|f(x)|\leq g(x)$ holds for any sufficiently large $x$. Therefore, we have
\[
\psi (x; q,a)\geq (1-e^{-aD})(1+o(1))\frac{x}{\varphi (q)}
\]
provided that $D>0$. Hence one can take $D_{PAP}=10c_{ZD}$, $C_{PAP}=1-e^{-aD_{PAP}}$.

Next, we evaluate $c_{ZD}$ above.  Put
\[
N^{*}(\alpha ,T,Q):=\sum _{q\leq Q}\sum _{\chi (\mathrm{mod}\; q)}^{\quad \quad *}N(\alpha ,T, \chi).
\]
By Theorem 1 of Jutila \cite{J}, for any fixed $\varepsilon >0$, 
\[
N^{*}(\alpha ,T,Q)\ll _{\varepsilon}(Q^{2}T)^{(2+\varepsilon)(1-\alpha)}
\]
holds uniformly for $\frac{4}{5}\leq \alpha \leq 1$, $T\geq 1$. Put $Q=T$. Then it follows that
\begin{equation}
\label{32}
\sum _{q\leq T}\sum _{\chi \; (\mathrm{mod}\; q)}^{\quad \quad *}N(\alpha ,T, \chi)\ll _{\varepsilon}T^{(6+\varepsilon)(1-\alpha)}.
\end{equation}
On the other hand, if $0\leq \alpha \leq \frac{4}{5}$, by using the trivial estimate 
\[
N (\alpha , T, \chi)\ll T\log T,
\]
we have 
\begin{equation}
\label{33}
\begin{aligned}
\sum _{q\leq T}\sum _{\chi \; (\mathrm{mod}\; q)}^{\quad \quad *}N(\alpha ,T, \chi) & \ll \sum _{q\leq T}qT\log T \ll T^{3}\log T \ll T^{(15+\varepsilon )(1-\alpha)}.
\end{aligned}
\end{equation}
By (\ref{32}) and (\ref{33}), 
\[
\sum _{q\leq T}\sum _{\chi \; (\mathrm{mod}\; q)}^{\quad \quad *}N(\alpha ,T, \chi)  \ll _{\varepsilon} T^{(15+\varepsilon)(1-\alpha)}
\]
holds uniformly for $0\leq \alpha \leq 1$. Therefore, one can take $c_{ZD}=16$. Consequently, one can take
\begin{equation}
\label{34}
C_{PAP}=1-e^{-2}, \quad D_{PAP}=160
\end{equation}
as the values of $C_{PAP}$ and $D_{PAP}$. 

%%%%%%%%%%%%%%%%%%%%%%%%%%%%%%%%%%%%%%%%%%%%%%%%%%%%%%%%%%%%%%%%%%%%%%%%%%%%%%%%%%%%%%%%%%%%%%%%%%%%%%%%%%%%
%%%%%%%%%%%%%%%%%%%%%%%%%%%%%%%%%%%%%%%%%%%%%%%%%%%%%%%%%%%%%%%%%%%%%%%%%%%%%%%%%%%%%%%%%%%%%%%%%%%%%%%%%%%%

\section{Proof of Proposition 3.1}
Proposition 3.1 can be derived from the following proposition.
\begin{prop}[An explicit version of Theorem 3 of \cite{FMT}]
Let $A\geq 1$ be an arbitrary real number and $x$ be a sufficiently large real number and $B_{0}$ be a positive integer. Put 
\[
y=c\frac{x\log x \log _{3}x}{\log _{2}x}, \quad c=\frac{\theta c_{I,J}}{12800 \log 5}
\]
and let $S(\vec{a})$, $S(\vec{n})$ be the subsets of $\mathbb{Z}$ defined in Section $2$. Then, there exist $A^{\prime}=A^{\prime}(x)$ and $c^{\prime}=c^{\prime}(x)$ satisfying 
\[
A^{\prime}\sim c^{\prime}A, \quad 1\leq c^{\prime} \leq 5
\]
and choices of vectors $\vec{\mathbf{a}}=(\mathbf{a}_{s}\; \mathrm{mod}\; s)_{s\in {\cal S}}$, $\vec{\mathbf{n}}=(\mathbf{n}_{p}\; \mathrm{mod}\; p)_{p\in {\cal P}}$ such that for arbitrarily fixed $0\leq \alpha <\beta \leq1$,
\begin{equation}
\label{Prop2}
\# ({\cal Q}\cap S(\vec{\mathbf{a}})\cap S(\vec{\mathbf{n}})\cap (\alpha y, \beta y])\sim A^{\prime}|\beta-\alpha|\frac{x}{\log x}
\end{equation}
holds with probability $1-o(1)$. The implied constant of $o(1)$ in the probability may depend on $A, \alpha ,\beta$, but that of (\ref{Prop2}) doesn't. 
\end{prop}
\noindent
({\it Proof of Proposition 6.1 $\Rightarrow$ Proposition 3.1}) Let $0<\varepsilon <1$. We decompose the interval $(0, 1]$ into $O(\varepsilon ^{-1})$ disjoint intervals $(\alpha _{i}, \beta _{i}]_{i\in I}$ of lengths between $\varepsilon /2$ and $\varepsilon$. We apply Proposition 6.1 with $(\alpha ,\beta)=(\alpha _{i}, \beta _{i})$ ($i\in I$) and $(\alpha ,\beta)=(0,1)$. Then there exist $A^{\prime}$, $c^{\prime}$ with
\[
A^{\prime}\sim c^{\prime}A, \quad 1\leq c^{\prime}\leq 5
\]
and some vectors of residue classes $\vec{a}=(a_{s}\; \textrm{mod}\; s)_{s\in {\cal S}}$, $\vec{n}=(n_{p}\; \textrm{mod}\; p)_{p\in {\cal P}}$ for which 
\begin{equation}
\label{35}
\# ({\cal Q}\cap S(\vec{a}) \cap S(\vec{n}))\sim A^{\prime}\frac{x}{\log x}
\end{equation}
and 
\begin{equation}
\label{36}
\# ({\cal Q}\cap S(\vec{a}) \cap S(\vec{n})\cap (\alpha _{i}y, \beta _{i}y])\leq A^{\prime}\varepsilon (1+o(1)) \frac{x}{\log x}
\end{equation}
hold. For any real numbers $0\leq \alpha <\beta \leq 1$, the interval $(\alpha y, \beta y]$ can be covered by at most $\lfloor \frac{2}{\varepsilon}\rfloor |\beta -\alpha |+1$ intervals $(\alpha _{i}y, \beta _{i}y]$. Hence by (\ref{36}) we have 
\begin{equation}
\label{37}
\begin{aligned}
\# ({\cal Q}\cap S(\vec{a}) \cap S(\vec{n}) \cap (\alpha y, \beta y])&\leq A^{\prime}\varepsilon \left( \left\lfloor \frac{2}{\varepsilon} \right\rfloor |\beta -\alpha |+1 \right)(1+o(1))\frac{x}{\log x} \\
&\leq A^{\prime}(2|\beta -\alpha |+\varepsilon )(1+o(1))\frac{x}{\log x}.
\end{aligned}
\end{equation}
We extend $\vec{a}$ to $(a_{p})_{p\leq x}$ by
\[
a_{p}:=\begin{cases} n_{p} \quad &(p\in {\cal P}) \\
0  \quad &(p\not \in {\cal S}\cup {\cal P}),
\end{cases}
\]
and set 
\[
{\cal T}:=\{n \in [y]\backslash [x]\; | \; n \not \equiv a_{p}\; (\mathrm{mod}\; p), \forall p\leq x, p\neq B_{0} \}.
\]
Following the argument of Section 4 of \cite{FMT}, we see that the set ${\cal T}$ only differs from ${\cal Q}\cap S(\vec{a}) \cap S(\vec{n})$ by a set ${\cal R}$ consisting of $z$-smooth numbers in $[y]$ multiplied by powers of $B_{0}$, and that ${\cal R}=o(x/\log x)$.  Consequently we obtain Proposition 3.1. \hspace{5.6cm} $\Box$ \\
Hence to prove the main theorem, it suffices to prove Proposition 6.1. This proposition can be derived from the following two propositions.
\begin{prop}[An explicit version of Theorem 5 of \cite{FMT}]
Let $x$ be a sufficiently large real number and put
\[
y=c\frac{x\log x \log _{3}x}{\log _{2}x}, \quad c=\frac{\theta c_{I,J}}{12800 \log 5}.
\]
Then, there exist a quantity $C=C(x)$ satisfying 
\begin{equation}
\label{38}
\frac{\theta c_{I,J}}{9600c}(1+o(1))\leq C \leq \frac{\theta c_{I,J}}{4800c}(1+o(1))
\end{equation}
and choices of random vectors $\vec{\mathbf{a}}=(\mathbf{a}_{s}\; \mathrm{mod}\; s)_{s\in {\cal S}}$, $\vec{\mathbf{n}}=(\mathbf{n}_{p} \; \mathrm{mod}\; p)_{p\in {\cal P}}$, for which the following three conditions hold. \\
$\cdot$ For any vector $\vec{a}$ in the essential range of $\vec{\mathbf{a}}$, 
\begin{equation}
\label{39}
\mathbb{P}(q\equiv \mathbf{n}_{p}\; (\mathrm{mod}\; p) \; | \; \vec{\mathbf{a}}=\vec{a})\leq x^{-\frac{1}{2}-\frac{1}{10}}
\end{equation}
holds uniformly for all $p\in {\cal P}$. \\
$\cdot$ For any fixed $0\leq \alpha <\beta \leq 1$, 
\begin{equation}
\label{40}
\# ({\cal Q}\cap S(\vec{\mathbf{a}})\cap (\alpha y, \beta y])\sim 80c |\beta -\alpha| \frac{x}{\log x}\log _{2}x
\end{equation}
holds with probability $1-o(1)$. \\
$\cdot$ We say that a vector $\vec{a}$ in the essential range of $\vec{\mathbf{a}}$ is “good” if 
\begin{equation}
\label{41}
\sum _{p\in {\cal P}}\mathbb{P}(q\equiv \mathbf{n}_{p}\; (\mathrm{mod}\; p) \; | \; \vec{\mathbf{a}}=\vec{a})=C+O_{\leq}\left( \frac{1}{(\log _{2}x)^{2}} \right)
\end{equation}
holds for any $q\in {\cal Q}\cap S(\vec{\mathbf{a}})$ with at most $x/\log x\log _{2}x$ exceptionals. Then $\vec{\mathbf{a}}$ is “good” with probability $1-o(1)$. 
\end{prop}

\begin{prop}[\cite{FMT}, Theorem 4]
Let $x$ be a sufficiently large real number, and ${\cal P}^{\prime}, {\cal Q}^{\prime}$  sets of primes in $(\frac{x}{2}, x]$, $(x, \log x]$ with $\# {\cal Q}^{\prime}>(\log _{2}x)^{3}$, respectively. For each $p\in {\cal P}^{\prime}$, let $\mathbf{e}_{p}$ be a random subset of ${\cal Q}^{\prime}$ satisfying 
\begin{equation}
\label{42}
\# \mathbf{e}_{p}\leq r=O\left( \frac{\log x \log _{3}x}{\log _{2}^{2}x} \right) \quad (p\in {\cal P}^{\prime}).
\end{equation}
Assume the following two conditions. \\
$\cdot$ For any $p \in {\cal P}^{\prime}$, $q\in {\cal Q}^{\prime}$, 
\begin{equation}
\label{43}
\mathbb{P}(q\in \mathbf{e}_{p})\leq x^{-\frac{1}{2}-\frac{1}{10}}.
\end{equation}
$\cdot$ For any $q\in {\cal Q}^{\prime}$ with at most $\frac{1}{(\log _{2}x)^{2}}\# {\cal Q}^{\prime}$ exceptions, there exists some quantity $C$ (which is independent of $q$) with 
\begin{equation}
\label{44}
\frac{5}{4}\log 5 \leq C \ll 1
\end{equation}
such that 
\begin{equation}
\label{45}
\sum _{p\in {\cal P}^{\prime}}\mathbb{P}(q\in \mathbf{e}_{p})=C+O_{\leq}\left( \frac{1}{(\log _{2}x)^{2}} \right)
\end{equation}
holds. \\
Then, for any positive integer $m$ with 
\begin{equation}
\label{46}
m\leq \frac{\log _{3}x}{\log 5},
\end{equation}
there exist  random subsets $\mathbf{e}_{p}^{\prime}\subset {\cal Q}^{\prime}$ for each $p \in {\cal P}$ such that 
\[
\# \{q\in {\cal Q}^{\prime} \; | \; q\not \in \mathbf{e}_{p}^{\prime}, \forall p \in {\cal P}^{\prime} \}\sim 5^{-m}\# {\cal Q}^{\prime}
\]
holds with probability $1-o(1)$. More generally, for any subset ${\cal Q}^{\prime \prime}\subset {\cal Q}^{\prime}$ with $\# {\cal  Q}^{\prime \prime}\geq (\#{\cal Q}^{\prime})/\sqrt{\log _{2}x}$, 
\[
\# \{q\in {\cal Q}^{\prime \prime} \; | \; q\not \in \mathbf{e}_{p}^{\prime}, \forall p \in {\cal P}^{\prime} \}\sim 5^{-m}\# {\cal Q}^{\prime \prime}
\]
holds with probability $1-o(1)$. 
\end{prop}

\quad \\
\noindent
({\it Proof of Proposition 6.2 $\land$ Proposition 6.3 $\Rightarrow$ Proposition 6.1 })
Let $m$ be a positive integer with 
\begin{equation}
\label{47}
m\leq \frac{\log _{3}x}{\log 5},
\end{equation}
and put $A^{\prime}=5^{-m}80c\log _{2}x$. If 
\[
m=\left\lfloor  \frac{\log \left( \frac{80c\log _{2}x}{A} \right)}{\log 5}  \right\rfloor ,
\]
then $A\leq A^{\prime}\leq 5A$ holds. For this $m$ satisfy (\ref{47}), we need 
\[
1\ll \frac{80c}{A}\leq 1.
\]
Write
\[
c=\frac{A}{80k}\varepsilon _{0}, \quad \varepsilon _{0}>0.
\]
Then by (\ref{38}), it follows that
\[
\frac{k\theta c_{I,J}}{120A\varepsilon _{0}}(1+o(1)) \leq C \leq \frac{k\theta c_{I,J}}{60A\varepsilon _{0}}(1+o(1)).
\]
This $C$ must satisfy the condition (\ref{44}) in Proposition 6.3. Hence the $\varepsilon _{0}$ must satisfy 
\[
1\ll \varepsilon _{0}\leq \frac{k\theta c_{I,J}}{150A\log 5}(1+o(1)).
\]
We take $\varepsilon _{0}=\frac{k\theta c_{I,J}}{160A\log 5}$. Then 
\[
c=\frac{k\theta c_{I,J}}{12800 \log 5}.
\]
Following the argument in Section 5 of \cite{FMT}, we see that the statement of Proposition 6.1 is valid with this $c$.  \hspace {12cm} $\Box$

With the above discussion, all that remains is to prove Proposition 6.2. In practice, however, we only need to check how the remaining coefficients are determined following the arguments in \cite{FGKMT} and \cite{FMT}. \\
\quad \\
\noindent
({\it Proof of Proposition 6.2}) Proposition 6.2 can be obtained by applying the sieve of Maynard (Theorem 6 of \cite{FMT}.  For the proof, see Section 6 of \cite{FGKMT}).  In particular, the constant $c_{0}$ in (6.1) of \cite{FMT} can be chosen as $c_{0}=\frac{1}{5}$, following the argument of Maynard \cite{May}. The quantity $C$ in Theorem 5 of \cite{FMT} is computed in the following way. Put
\[
\sigma =\prod _{s\in {\cal S}}\left(1-\frac{1}{s} \right), \quad y=c\frac{x\log x \log _{3}x}{\log _{2}x}.
\]
By Mertens formula, we have
\[
\sigma y =\left(1+O\left( \frac{1}{\log _{2}^{10}x} \right) \right)80cx\log _{2}x.
\]
Following (8.3) of \cite{FGKMT}, we take
\[
u=\frac{\varphi (B)}{B}\frac{\log R}{\log x}\frac{rJ_{r}}{2I_{r}},
\]
where $B$ is either $1$ or prime, so
\[
\frac{1}{2}\leq \frac{\varphi (B)}{B}\leq 1.
\]
Following Theorem 6 of \cite{FGKMT}, we take $R=x^{\frac{\theta}{3}}$, where $\theta >0$ is a constant which satisfies the condition in Definition 2 of \cite{FGKMT}. (As mentioned above, one can take $\theta =\frac{1}{3}$.) Suppose
\[
J_{r}=c_{I,J}(1+o(1))\frac{\log r}{r}I_{r}
\]
holds as $r\to \infty$. Then
\[
\frac{\theta}{12}c_{I,J}(1+o(1))\log r \leq u \leq \frac{\theta}{6}c_{I,J}(1+o(1))\log r.
\]
The parameter $r$ above can be chosen as 
\[
r=\left\lfloor \log ^{c_{0}}x \right\rfloor = \left\lfloor \log ^{\frac{1}{5}}x \right\rfloor .
\]
Hence
\[
\frac{\theta}{60}c_{I,J}(1+o(1))\log _{2}x \leq u \leq \frac{\theta}{30}c_{I,J}(1+o(1))\log _{2}x.
\]
Finally, following \cite{FMT}, one can take
\[
C=\frac{u}{\sigma}\frac{x}{2y}.
\]
Hence due to the above estimate of $u$, the statement of Proposition 6.1 holds if the parameter satisfies the condition  (\ref{38}).  \hspace{9.7cm}  $\Box$

%%%%%%%%%%%%%%%%%%%%%%%%%%%%%%%%%%%%%%%%%%%%%%%%%%%%%%%%%%%%%%%%%%%%%%%%%%%%%%%%%%%%%%%%%%%%%%%%%%%%%%%%%%%%
%%%%%%%%%%%%%%%%%%%%%%%%%%%%%%%%%%%%%%%%%%%%%%%%%%%%%%%%%%%%%%%%%%%%%%%%%%%%%%%%%%%%%%%%%%%%%%%%%%%%%%%%%%%%

\section{Acknowledgements}
This work is partially supported by the JSPS, KAKENHI Grant Number 24K06697. The author would like to thank Professor Terence Tao for telling him an idea to obtain a lower bound (\ref{Tao})  in his blog. 
%%%%%%%%%%%%%%%%%%%%%%%%%%%%%%%%%%%%%%%%%%%%%%%%%%%%%%%%%%%%%%%%%%%%%%%%%%%%%%%%%%%%%%%%%%%%%%%%%%%%%%%%%%%%
%%%%%%%%%%%%%%%%%%%%%%%%%%%%%%%%%%%%%%%%%%%%%%%%%%%%%%%%%%%%%%%%%%%%%%%%%%%%%%%%%%%%%%%%%%%%%%%%%%%%%%%%%%%%

\noindent
Kanto Gakuin University, \\
Kanazawa, Yokohama,\\
Kanagawa, Japan\\
E-mail address: sono@kanto-gakuin.ac.jp

\end{document}